\documentclass[reqno,11pt]{amsart}
\usepackage{amsthm,amsfonts,amssymb,euscript,mathrsfs,graphics,xcolor,amsmath,latexsym,marginnote,hyperref}
\usepackage[latin1]{inputenc}
\usepackage{amsaddr}
\usepackage{overpic}
\usepackage{graphicx}

\theoremstyle{plain}

\usepackage{times}
\numberwithin{equation}{section}

\numberwithin{equation}{section}

\theoremstyle{plain}
\newtheorem{teor}{Theorem}[section]
\newtheorem{ese}[teor]{Example}
\newtheorem{prop}[teor]{Proposition}
\newtheorem{lem}[teor]{Lemma}
\newtheorem{cor}[teor]{Corollary}

\newcommand{\bdm}{\begin{displaymath}}
\newcommand{\edm}{\end{displaymath}}
\newcommand{\bpb}{\begin{prob}}
\newcommand{\epb}{\end{prob}}

\newcommand{\beq}{\begin{equation}}
\newcommand{\eeq}{\end{equation}}
\newcommand{\bem}{\begin{multline}}
\newcommand{\eem}{\end{multline}}
\newcommand{\bes}{\begin{ese}}
\newcommand{\ees}{\end{ese}}
\newcommand{\bde}{\begin{defi}}
\newcommand{\ede}{\end{defi}}
\newcommand{\bpr}{\begin{prop}}
\newcommand{\epr}{\end{prop}}
\newcommand{\ble}{\begin{lem}}
\newcommand{\ele}{\end{lem}}
\newcommand{\bte}{\begin{teor}}
\newcommand{\ete}{\end{teor}}
\newcommand{\bco}{\begin{cor}}
\newcommand{\eco}{\end{cor}}

\theoremstyle{definition}
\newtheorem{defi}[teor]{Definition}
\newtheorem{remark}[teor]{Remark}

\def\Z{\mathbb{Z}}

\newcommand{\eps}{\varepsilon}

\def\e{\varepsilon}
\def\N{{\mathbb N}}
\renewcommand{\to}{\rightarrow}

\newcommand{\T}{\mathbb{T}}

\newcommand{\NN}{{\mathcal N}}

\newcommand{\calP}{{\mathcal P}}

\newcommand{\TT}{{\mathcal T}}

\newcommand{\al}{\alpha}

\newcommand{{\resonance}}{relevant self-energy cluster }

\newcommand{\de}{\delta}
\newcommand{\ZZ}{\mathbb{Z}}
\newcommand{\OO}{\mathcal{O}}

\newcommand{\odd}{\mathrm{odd}}

\makeatletter

\def\l@subsection{\@tocline{2}{0pt}{2.5pc}{5pc}{}}
\def\l@subsubsection{\@tocline{3}{0pt}{4.5pc}{5pc}{}}

\renewcommand\tocchapter[3]{%
  \indentlabel{\@ifnotempty{#2}{\ignorespaces#2.\quad}}#3%
}
\newcommand\@dotsep{4.5}
\def\@tocline#1#2#3#4#5#6#7{\relax
  \ifnum #1>\c@tocdepth 
  \else
    \par \addpenalty\@secpenalty\addvspace{#2}%
    \begingroup \hyphenpenalty\@M
    \@ifempty{#4}{%
      \@tempdima\csname r@tocindent\number#1\endcsname\relax
    }{%
      \@tempdima#4\relax
    }%
    \parindent\z@ \leftskip#3\relax \advance\leftskip\@tempdima\relax
    \rightskip\@pnumwidth plus1em \parfillskip-\@pnumwidth
    #5\leavevmode\hskip-\@tempdima{#6}\nobreak
    \leaders\hbox{$\m@th\mkern \@dotsep mu\hbox{.}\mkern \@dotsep mu$}\hfill
    \nobreak
    \hbox to\@pnumwidth{\@tocpagenum{#7}}\par
    \nobreak
    \endgroup
  \fi}
\def\l@subsection{\@tocline{2}{0pt}{2.5pc}{5pc}{}}

\begin{document}

\begin{center}
{\small Scientific chapter: 18. Partial Differential Equations }
\end{center}
\bigskip

\title{ Chaotic resonant dynamics and exchanges of energy in Hamiltonian PDEs }

\date{}


\author{Filippo Giuliani}
\address{UPC, Barcelona, \\e-mail: filippo.giuliani@upc.edu}

\author{Marcel Guardia}
\address{UPC, Barcelona, \\e-mail: marcel.guardia@upc.edu}

\author{Pau Martin}
\address{UPC, Barcelona, \\e-mail: p.martin@upc.edu}

\author{Stefano Pasquali}
\address{Matematikcentrum, Lunds Universitet, Lund, \\e-mail: stefano.pasquali@math.lu.se}


 \begin{abstract}
 The aim of this note is to present the recent results in \cite{GGMP} where we provide the existence of solutions of some nonlinear resonant PDEs on $\mathbb{T}^2$ exchanging energy among Fourier modes in a ``chaotic-like'' way.  We say that a transition of energy is ``chaotic-like'' if either the
choice of activated modes or the time spent in each transfer
can be chosen randomly. We consider the nonlinear cubic Wave, the Hartree and the nonlinear cubic
Beam equations.
The key point of the construction of the special solutions is the existence of heteroclinic connections between invariant objects and the
construction of symbolic dynamics (a Smale horseshoe) for the Birkhoff
Normal Form of those equations.
\end{abstract}

\maketitle
\section{Introduction}

A fundamental question in nonlinear Hamiltonian Partial Differential Equations (PDEs) on compact manifolds is to understand how solutions can exchange energy among Fourier modes as time evolves.

In the last decade there has been a lot of activity in building exchange of energy behaviors in different Hamiltonian PDEs almost exclusively  for the nonlinear Schr\"odinger equation. They can be classified into two groups: the first one are the so-called \emph{beating solutions} \cite{GrebertV11,GT,GPT,HT, HausP17}, namely orbits that are essentially supported on a finite numbers of modes and whose energy oscillates between those modes in a certain time range; the other group are those addressing the problem of \emph{transfer of energy}, namely constructing orbits whose energy is transferred to increasingly higher modes as time evolves \cite{Bourgain96,Kuksin96,Kuksin97b,CKSTT,Hani12,Guardia14,GuardiaK12,HaniPTV15, HausProcesi,GuardiaHP16, Pocovnicu11, Pocovnicu12, Maspero18g, Delort10, GerardG10,  GerardG11}.

Most of  these results rely on analyzing the first order Birkhoff normal form of the Hamiltonian PDEs and building invariant objects for such models. Note that these first order Birkhoff normal forms are typically  non-integrable Hamiltonian systems (at least in dimension greater or equal than 2). Nevertheless, restricted to suitably chosen invariant subspaces those models are integrable (they have ``enough'' first integrals in involution). The analysis of the dynamics at (or near) these invariant subspaces allows to construct unstable motions and exchange of energy solutions.

The purpose of this note is to present the recent results in \cite{GGMP}
where we consider three different PDEs, a nonlinear Wave equation, a nonlinear Beam equation and the Hartree equation (see \eqref{Wave}, \eqref{Beam} and \eqref{Hartree} below) and we show the existence of solutions that display exchange of energy behaviors in a chaotic fashion (up to a certain time scale). They are thus rather  different from the previously constructed beating solutions results \cite{GT,GPT,HT}, whose exchange of  energy is periodic in time.
On the contrary,  the beating solutions we construct undergo oscillations that can be ``randomly'' chosen (see Section \ref{sec:NLWMain} for the precise statements).
This ``random'' choice is obtained by exploiting the non-integrability and chaoticity (symbolic dynamics) of its Birkhoff normal form.

\subsection{Main results}\label{sec:NLWMain}
Consider the completely resonant cubic nonlinear Wave and Beam equations on the $2$-dimensional torus
\begin{align}
u_{tt}-\Delta u+u^3&=0 \qquad u=u(t, x), \quad t\in\mathbb{R}, \quad x\in\mathbb{T}^2\label{Wave}\\
u_{tt}+\Delta^2 u+u^3&=0\qquad u=u(t, x), \quad t\in\mathbb{R}, \quad x\in\mathbb{T}^2.\label{Beam}
\end{align}
We prove the existence of special beating solutions for such PDEs
which are (essentially) Fourier supported on a finite set of $4$-tuple resonant modes
 \begin{equation}\label{def:lambdaset}
 \Lambda:=\{n^{(r)}_j\}^{r=1, \dots, N}_{j=1, \dots, 4}\subset\Z^2,\quad
N\geq 2,
 \end{equation}
in the sense that
 \[
 u(t, x)=\sum_{j\in\Lambda} a_j(t)\,e^{\mathrm{i} j\cdot x}+R(t, x)
 \]
 where $R(t, x)$ is small in some Sobolev norm. The transfers of energy between modes in $\Lambda$ are \emph{chaotic-like}, in the following sense. Either
\begin{itemize}
\item[(a)] one can prescribe a finite sequence of times $t_1, \dots, t_n$ and find a solution that exists for \emph{long but finite time} exhibiting transfers of energy among the modes in $\Lambda$ at the prescribed times $t_1, \dots, t_n$
\end{itemize}
 or
 \begin{itemize}
 \item[(b)] one can prescribe any sequence of resonant tuples
$\{n^{(r_h)}_{j}\}_{h=1, \dots, k}\subseteq \Lambda$ and find a  solution and a
sequence of times $t_1, \dots, t_k$ such that at time $t_n$ all modes are
"switched off" (modulus of the modes almost constant) whereas the
modes $(n^{(r_n)}_{1},n^{(r_n)}_{2},n^{(r_n)}_{3},n^{(r_n)}_{4})$ are ``switched
on", in the sense that they exchange energy between them.
\end{itemize}
We look for these solutions in the subspace
\begin{align*}
\mathcal{U}_{\odd} = \left\{ u=\sum_{j\in\mathbb{Z}_{\odd}^2} u_j e^{\mathrm{i} j \cdot x} \right  \}, \quad
\mathbb{Z}^2_{\odd}= \left\{ (j^{1},
j^{2})\in\mathbb{Z}^2 : \,
j^{1}\text{ odd},j^{2}\mbox{ even}\right \},
\end{align*}
which is invariant under the flow of  equations \eqref{Wave}, \eqref{Beam} (see \cite{Procesi2010ANF}).
The origin of such subspace is an elliptic fixed point, and solutions of the variational equation
\begin{equation}\label{def:variational}
\ddot{u}_j+\omega^2(j) u_j=0, \quad j\in\mathbb{Z}^2_\odd,
\end{equation}
where $\omega(j)=|j|$ (for the Wave equation \eqref{Wave}) and
$\omega(j)=|j|^2$ (for the Beam equation \eqref{Beam}),
 are superposition of decoupled harmonic oscillators. Hence all solutions are periodic/quasi-periodic/almost-periodic in time and,  in particular, there is no transfer of energy between the linear modes when time evolves.
This implies that the existence of beating solutions depends on the presence of the nonlinearities.
In order to describe the nonlinear effects in a neighborhood of an elliptic equilibrium we perform a Birkhoff normal form analysis: namely, we construct changes of coordinates that transform the Hamiltonian of the equations \eqref{Wave}, \eqref{Beam} into a Hamiltonian of the form
\begin{equation}\label{def:BNFAbs}
K=K^{(2)}+K^{(4)}+\mathcal{R},
\end{equation}
where $K^{(i)}$ are homogenous terms of degree $i$ and $\mathcal{R}$ is a function that can be considered as a small perturbation. Then, one can consider the truncated system
\begin{equation}\label{def:HamN}
\NN:=K^{(2)}+K^{(4)},
\end{equation}
called \emph{normal form},  as a model which describes the effective dynamics of equations \eqref{Wave}, \eqref{Beam} up to a certain time scale. 
It is well known that the existence of such changes of coordinates cannot be always guaranteed because of the presence of small divisor problems and/or derivatives in the nonlinear terms. To overcome these problems we adopt the strategy of performing a \emph{weak version} of the Birkhoff normal form (see Section \ref{sec:Proof}) which is well established, for instance, in the KAM theory for quasi-linear resonant PDEs \cite{GiulianiKdV}, \cite{FGP}.

The normal form Hamiltonian $\NN$ possesses many finite-dimensional, symplectic, invariant subspaces of the form
\begin{equation}\label{def:Vlambda}
V_\Lambda:=\{ u_j=0\,\, \,\,\forall j\notin \Lambda \},
\end{equation}
 where $\Lambda\subset\mathbb{Z}^2_\odd$ is a finite set (suitably chosen).\\
 In the following theorem we state two results concerning the dynamics of the normal form Hamiltonian $\mathcal{N}$, which are fundamental ingredients for the proof of the main Theorems \ref{TeoWaveBeam} and \ref{thm:traveling_periodic_beating}.
 
To state this theorem, let us first introduce the Bernoulli shift. Consider $\Sigma=\N^\ZZ$, the space of sequences of natural numbers, with the usual topology taking as neighborhood basis of $\omega^*=\{\omega^*_j\}_{j\in\mathbb{Z}}$ the sets
\[
 U_k=\left\{\omega\in\Sigma: \omega_j=\omega_j^*\quad \text{for}\quad |j|<k\right\}.
\]
Then,   the Bernouilli shift is a homeomorphism defined as
\begin{equation}\label{def:shift}
 \sigma:\Sigma\to\Sigma,\qquad (\sigma\omega)_k=\omega_{k+1}.
\end{equation}
The map $\sigma$ is one of the paradigmatic examples of chaotic dynamics and encodes the dynamics of the classical Smale Horseshoe (of infinite symbols), see~\cite{Moser01}. In particular it has dense orbits, its periodic orbits form a dense set in $\Sigma$ and  it has positive topological entropy.
\begin{teor}\label{thm:geometric}
Let $N\geq 2$. There exist sets $\Lambda\subset\mathbb{Z}^2_\odd$ of cardinality $4 N$ such that $V_\Lambda$ is invariant by the dynamics of $\NN$ and the following holds.
\begin{itemize}
\item[(i)] Let $N=2$. Then, the flow $\Phi_t$ associated to $\NN$ in $V_\Lambda$ has the following property. There exists a section $\Pi$ transverse to the flow $\Phi_t$ such that the induced Poincar\'e map
\[
 \calP:\mathcal{U}=\mathring{\mathcal{U}}\subset\Pi\to\Pi
\]
has an invariant set $X\subset \mathcal{U}$ which is homeomorphic to $\Sigma\times\T^5$.
Moreover, the dynamics of $\calP:X\to X$ is topologically conjugated to the map
\[
 \widetilde\calP:\Sigma\times\T^5\to \Sigma\times\T^5, \qquad  \widetilde\calP(\omega, \theta)=(\sigma\omega, \theta+f(\omega))
\]
where $\sigma$ is the Bernoulli shift \eqref{def:shift} and $f:\Sigma\to\T^5$ is a continuous function (see Remark \eqref{rmk:factor} below).

Namely, $\calP$ has a Smale horseshoe of infinite symbols as a factor.
\item[(ii)] There exist $N$ partially hyperbolic  $2(N+1)$-dimensional tori, $\mathbb{T}_1, \dots, \mathbb{T}_N$, invariant for the restriction of the normal form Hamiltonian $\mathcal{N}$ at the subspace $V_\Lambda$, which have the following property. Take arbitrarily small neighborhoods $V_i$ of $\mathbb{T}_i$ and  any sequence $\{p_i\}_{i\geq 1}\subset \{ 1, \dots, N \}^\N$. Then, there exists an orbit $v(t)$ of $\mathcal{N}$ and a sequence of times $\{t_i\}_{i\geq 1}$ such that
\[
v(t_i)\in V_{p_i}.
\]
\end{itemize}
\end{teor}

\begin{remark}\label{rmk:factor}
The normal form Hamiltonian restricted to $V_{\Lambda}$ has $5$ constants of motion which are in involution and linearly independent in a neighborhood of the invariant set $X$. This implies the existence of a set of coordinates that puts these first integrals as actions. The angles $\theta$ appearing in Theorem \ref{thm:geometric} are just their conjugate variables. In particular the angles are cyclic, i.e. the dynamics of the angles on the section $\Pi$ is given by a translation $f$ which does not depend on $\theta$ but only on $\omega$.
\end{remark}

\begin{remark}\label{rmk:Lambda}
The set $\Lambda\subset\ZZ_{\odd}^2$ is the union of $N$ resonant tuples (with
certain
properties). The ``shape'' of the resonant tuples
$n_1,n_2,n_3,n_4\in\mathbb{Z}^2_{\odd}$ are
different for the Beam and Wave Equations. For the Beam equation, as for the
cubic nonlinear Schr\"odinger equation, they are rectangles with vertices in
$\ZZ^2_{\odd}$, since they must satisfy
\[
 n_1-n_2+n_3-n_4=0,\quad | n_1|^2-|n_2|^2+|n_3|^2-|n_4|^2=0.
\]
For the Wave equation they satisfy
\[
 n_1-n_2+n_3-n_4=0,\quad | n_1|-|n_2|+|n_3|-|n_4|=0.
\]
Those tuples form a parallelogram inscribed on an ellipse with foci at $F_1=0$
and $F_2=n_1+n_3$ and semi-major axis $a=(|n_1|+|n_3|)/2$.

Actually there is a ``\emph{large}'' choice for the sets $\Lambda\subset\ZZ^2_{\odd}$
for which Theorem \ref{thm:geometric} (and also Theorems \ref{TeoWaveBeam} and
\ref{thm:traveling_periodic_beating} below) is satisfied. Indeed,  Theorem
\ref{thm:geometric} relies on  the existence of a transverse intersection between certain invariant manifolds. This transversality is proven by  perturbative methods and, therefore, we need the restriction of $\mathcal{N}$ on $V_{\Lambda}$ (see \eqref{def:HamN}, \eqref{def:Vlambda}) to be close to integrable. This relies on choosing appropriate sets $\Lambda$.
The precise statement goes as follows. Fix $\eps>0$ (which will measure the
closeness to integrability). Then, for any $R\gg 1$,  one can  choose the
resonant tuples in the set $\Lambda$ generically in the annulus
 \[
  R(1-\eps)\leq |n|\leq R (1+\eps).
 \]
Generically means that one has to exclude the zero set of a finite number of
algebraic varieties (and the number of those is independent of $\eps$ and $R$).
\end{remark}

Items $(a)$ and $(b)$ above are consequence  of items $(i)$ and
$(ii)$ in Theorem \ref{thm:geometric}, respectively. Let us make some remark on the type of
dynamics for the normal form Hamiltonian $\mathcal{N}$.
\begin{itemize}
 \item Item (i) of Theorem \ref{thm:geometric} gives the existence of an invariant
set for $\mathcal{N}$ (see \eqref{def:HamN}) which possesses chaotic dynamics.
Such chaotic dynamics is obtained through the classical Smale horseshoe dynamics
for a suitable Poincar\'e map. This invariant set is constructed in the
neighborhood of homoclinic points to an invariant tori orbit (which, after a suitable symplectic reduction, becomes a
periodic orbit with large period $\mathtt{T}\gg 1$). The (infinite) symbols
codify the closeness to the invariant manifolds of the periodic orbit, the larger the symbol is, the longer the return time to the section
$\Pi$ is. In particular, one can construct orbits which take longer and longer
time to return to $\Pi$ for higher iterates. 

 Even if the theorem, as stated, gives the existence of one invariant set, one
actually can construct  a Smale horseshoe at each energy level.

 \item Item (ii) of Theorem \ref{thm:geometric} gives orbits which visit
(possibly infinitely many times) a given set of invariant tori in any prescribed
order. The construction of such orbits follows the classical strategy of Arnold
Diffusion \cite{Arnold64}. That is, it is a consequence of the existence of a chain
of invariant tori (again periodic orbits in a suitable symplectic reduction)
connected by transverse heteroclinic connections, plus a classical shadowing
argument (Lambda lemma, see for instance \cite{FontichM98}).


 As for item (i) one can obtain the explained behavior at each energy level.
Indeed, the invariant tori come in families parameterized by the energy level
and therefore one can obtain this shadowing behavior at each energy level as
well.

\end{itemize}

Note that the knowledge of the orbits obtained in Theorem \ref{thm:geometric} is
\emph{for all time}. If one adds the error terms dropped from the original equation,
that is $\mathcal{R}$ in \eqref{def:BNFAbs}, one can obtain orbits for equations
\eqref{Wave}, \eqref{Beam} which follow the orbits of  Theorem
\ref{thm:geometric} for some time scales. Next theorem gives solutions of
equations \eqref{Wave} and \eqref{Beam} which (approximately) behave as those
obtained in Item $(i)$ of Theorem \ref{thm:geometric}.

\begin{teor}\label{TeoWaveBeam}
Let $N=2$ and fix  $0<\varepsilon\ll 1$. Then for a large choice of sets
$\Lambda=\{n_i\}_{i=1}^8\subset\mathbb{Z}^2$ as in \eqref{def:lambdaset}
there exists $\mathtt{T}_0\gg 1$ such that for all $\mathtt{T}\geq \mathtt{T}_0$
there exists $M_0>0$ such that for all $M\geq M_0$ there exists
$\delta_0=\delta_0(M, \varepsilon, \mathtt{T})>0$ such that $\forall \delta\in
(0, \delta_0)$ the following holds.

Choose any   $k\geq  1$ and any sequence $\{ m_j\}_{j=1}^k$ such that $m_j\geq
M_0$ and $\sum_{j=1}^k m_j\leq M-k$.  Then, there exists a
solution $u(t, x)$ of \eqref{Wave}, \eqref{Beam} for $t\in [0, \delta^{-2} M
\mathtt{T}]$  of the form
\begin{equation*}
u(t, x)=\frac{\de}{\sqrt{2}}\sum_{i=1}^ 8 |n_i|^{-\kappa/2}\left(a_{n_i}(
t)\,e^{\mathrm{i} n_i \cdot x}+\overline{a_{n_i}}(t)\,e^{-\mathrm{i} n_i \cdot
x}\right)+R_1(t, x)
\end{equation*}
where $\kappa=1$ for the Wave equation \eqref{Wave} and $\kappa=2$ for the Beam
equation \eqref{Beam}, and  $\sup_{t\in [0, \delta^{-2} M \mathtt{T}]}\lVert R_1
\rVert_{H^s(\T^2)}\lesssim_s \delta^{3/2}$ for all $s\geq 0$. The first order
$\{ a_{n_i}\}_{i=1\dots 8}$
satisfies
\[
\begin{split}
\lvert a_{n_1}(t) \rvert^2&=\lvert a_{n_3} (t) \rvert^2=1- \lvert a_{n_2}(t)
\rvert^2=1-\lvert a_{n_4} (t) \rvert^2,\\
\lvert a_{n_5}(t) \rvert^2&=\lvert a_{n_7} (t) \rvert^2=1- \lvert a_{n_6}(t)
\rvert^2=1-\lvert a_{n_8} (t) \rvert^2,
\end{split}
\]
and has the following behavior.
\begin{itemize}
 \item \textbf{First resonant tuple (Periodic transfer of energy):} There exists
a $\mathtt{T}$-periodic function $Q(t)$, independent of $\delta$ and satisfying
$\min_{[0,\mathtt{T}]} Q(t) < \eps$ and $\max_{[0,\mathtt{T}]} |Q(t)| > 1-\eps$,
such that
 \[\lvert a_{n_1}(t) \rvert^2= Q (\delta^{2} t)+R_2(t)\qquad \text{with} \qquad
\sup_{t\in \mathbb{R}} \lvert R_2(t) \rvert\leq \varepsilon.\]
\item \textbf{Second resonant tuple (Chaotic-like transfer of energy):} There
exists a sequence of times $\{t_j\}_{j=0}^{k}$ satisfying $t_0=0$ and
\[
 t_{j+1}=t_j+\de^{-2}\mathtt{T}\left(m_j+\theta_j\right)\qquad \text{with}\qquad
\theta_j\in (0,1)
\]
such that
\[
 \lvert a_{n_5}( t_j) \rvert^2=\frac{1}{2}.
\]
Moreover, there exists another sequence $\{\bar t_j\}_{j=1\ldots k}$ satisfying
$t_j<\bar t_j<t_{j+1}$ such that,
\begin{equation}\label{eq:TheoRandomTimeOsc1}
\begin{split}
\lvert a_{n_5}(t) \rvert^2&>\frac{1}{2}\qquad \text{for}\qquad t\in (t_j,\bar
t_j)\\
\lvert a_{n_5}(t) \rvert^2&<\frac{1}{2}\qquad \text{for}\qquad t\in (\bar
t_j,t_{j+1})
\end{split}
\end{equation}
and
\begin{equation}\label{eq:TheoRandomTimeOsc2}
\sup_{t\in (t_j,\bar t_j)}\lvert a_{n_5}(t) \rvert^2\geq 1-\varepsilon\qquad
\text{ and }\qquad
\inf_{t\in (\bar t_j, t_{j+1})}\lvert a_{n_5}(t) \rvert^2\leq \varepsilon.
\end{equation}
\end{itemize}
\end{teor}

\begin{figure} 
\begin{overpic}[scale=0.9]{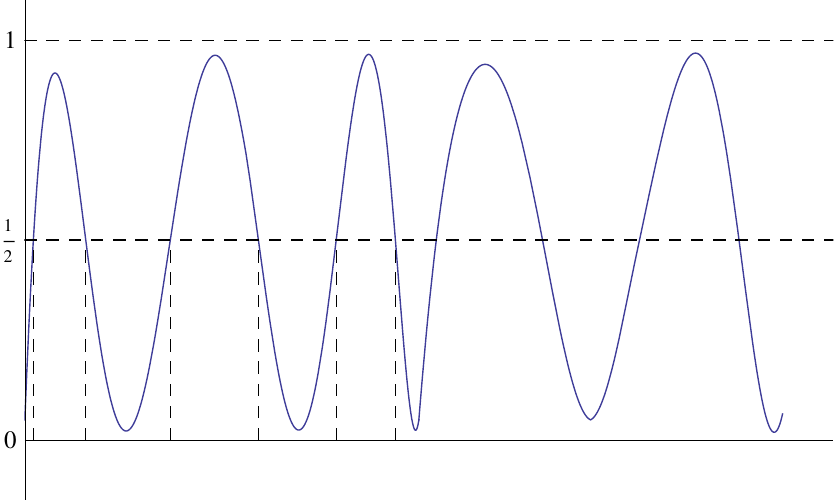}
\put(4,-.1){$t_1$}
\put(10,-.1){$\bar{t}_1$}
\put(18,-.1){$t_2$}
\put(30,-.1){$\bar{t}_2$}
\put(38,-.1){$t_3$}	
\put(46,-.1){$\bar{t}_3$}
\put(-15,60){$|a_{n_{5}}(t)|^2$}
\put(99,-.1){ $t$ }
\end{overpic}

\caption{An example of the evolution of the energy $|a_{n_5}(t)|^2$ as time evolves. The energy is a multi-bump like function. It assumes the value $1/2$ at the ``random" times $t=t_j$ and also at $t=\bar{t}_j$. The randomness in the $t_j$'s prescribes the separation of the bumps. The larger is the increment $t_{j+1}-t_j$, the more separated are the corresponding bumps. This shows that one can obtain very complicated energy transfer behaviors for the second resonant tuple.}
\label{fig:orbit}
\end{figure}
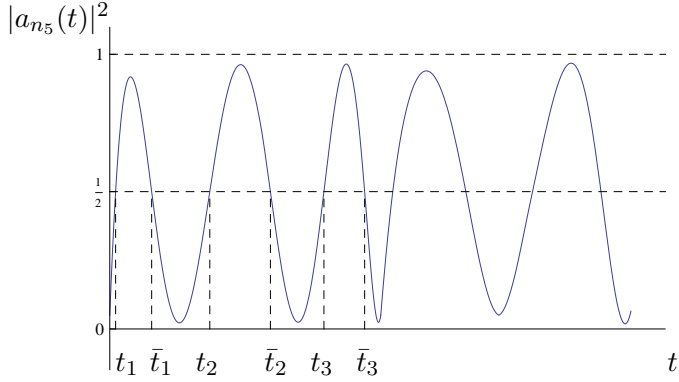

Note that the first order $\{\delta a_{n_i}\}_{i=1\ldots8}$ are the trajectories
obtained in Theorem \ref{thm:geometric}--(i) which belong to the horseshoe.

The first resonant tuple has a periodic beating behavior similar to \cite{GT}.
On the contrary, the behavior of the second resonant tuple is radically
different. The modulus of the modes $a_{n_i}$, $i=5,6,7,8$ ``oscillate'' from
being   $\OO(\eps)$ to  being $\OO(\eps)$-close to 1 (see Figure \ref{fig:orbit}). Moreover,  the time
intervals $\Delta t_j=t_j-t_{j-1}$ between times  $\{t_j\}$ when the
modes in the tuple have the same modulus, that
is
\[
\lvert a_{n_5}( t_j) \rvert^2=\lvert a_{n_6}( t_j) \rvert^2=\lvert a_{n_7}( t_j)
\rvert^2= \lvert a_{n_8}( t_j) \rvert^2=\frac{1}{2},
\]
(and the modulus of $a_{n_5}$ and $ a_{n_7}$ is increasing) can be chosen
randomly as any (large enough) integer multiple of $\mathtt{T}$ (plus a small
error).



Now we state the second main result of this paper, which gives solutions of
equations~\ref{Wave} and~\ref{Beam} which (approximately) behave as those
obtained in Item~$(ii)$ of Theorem~\ref{thm:geometric}.
\begin{teor}\label{thm:traveling_periodic_beating}
Let $N\geq 2$,  $k\gg 1$, $0<\eps\ll 1$.  Then for a large choice of a set
$\Lambda:=\{n^{(r)}_j\}^{r=1, \dots, N}_{j=1, \dots, 4} \subset\mathbb{Z}^2$ as
in \eqref{def:lambdaset} there exist $\delta_0>0$, $T>0$,  such that for any
$\de\in (0,\de_0)$ and any  sequence $\omega=(\omega_1, \dots, \omega_k),
\omega_i\in \{1, \dots, N  \}$, there exists
 a solution $u(t, x)$ of the \eqref{Wave}, \eqref{Beam}
 of the form
\begin{equation*}
u(t, x)=\frac{\de}{\sqrt{2}}\sum_{n\in \Lambda}
|n|^{-\kappa/2}\left(a_{n}(t)\,e^{\mathrm{i} n \cdot
x}+\overline{a_{n}}(t)\,e^{-\mathrm{i} n \cdot x}\right)+R_3(t, x), \quad t\in
[0, \delta^{-2} T]
\end{equation*}
where $\kappa=1$, for the Wave equation~\eqref{Wave}, or $\kappa=2$, for the Beam
equation~\eqref{Beam}, $ \sup_{t\in[0, \delta^{-2} T]} \| R_3(t, x)
\|_{H^s(\mathbb{T}^2)}\lesssim_s \delta^{3/2}$ for all $s\geq 0$,
and the  first order $\{ a_{n}\}_{n\in\Lambda}$,
has the following behavior.

There exist  some $\al_p,\beta_p$ satisfying
\[
\al_p<\beta_p<\al_{p+1} \qquad \text{ and }\qquad \beta_{p}-\alpha_p\gtrsim
|\ln\eps|, \quad p=1, \dots, k
\]
such that, if one splits the time interval as $[0, \delta^{-2} T]=I_1\cup J_{1,
2}\cup I_2\cup J_{2, 3}\cup\dots \cup J_{k-1, k }\cup I_k$
with
\[
I_p=[\delta^{-2}\al_p, \delta^{-2}\beta_p],\quad
J_{p,p+1}=[\delta^{-2}\beta_p,\delta^{-2}\al_{p+1}],
\]
such that $\{a_{n}\}_{n\in\Lambda}$ satisfies:
\begin{itemize}
\item in the \textbf{beating-time} intervals $I_p$, there exists $t_p>0$ such
that
\begin{align*}
\sup_{t\in I_p} \Big| |a_{n^{(\omega_p)}_1}(t)|^2-Q(\delta^2 t-t_p)
\Big|&\le\eps \\
\sup_{t\in I_p}  |a_{n^{(r)}_1}(t)|^2
& \leq \varepsilon \qquad &\text{ for }&\qquad r= 1, \dots, N, \quad r\neq \omega_p,
\end{align*}
where $Q(t)$ is the periodic function given by Theorem \ref{TeoWaveBeam},
\item in the \textbf{transition-time} intervals $J_{p, p+1}$,
\begin{align*}
\sup_{t\in J_{p, p+1}} |a_{n^{(r)}_1}(t)|^2&\geq 1-\eps \qquad \quad&\text{ for
}&\qquad r= 1, \dots, N,  \quad r\neq \omega_p,\\
\end{align*}
\end{itemize}
and $|a_{n^{(r)}_1}(t)|^2=|a_{n^{(r)}_3}(t)|^2$ ,
$|a_{n^{(r)}_j}(t)|^2=1-|a_{n^{(r)}_{1}}(t)|^2$ with $j=2, 4$.
\end{teor}

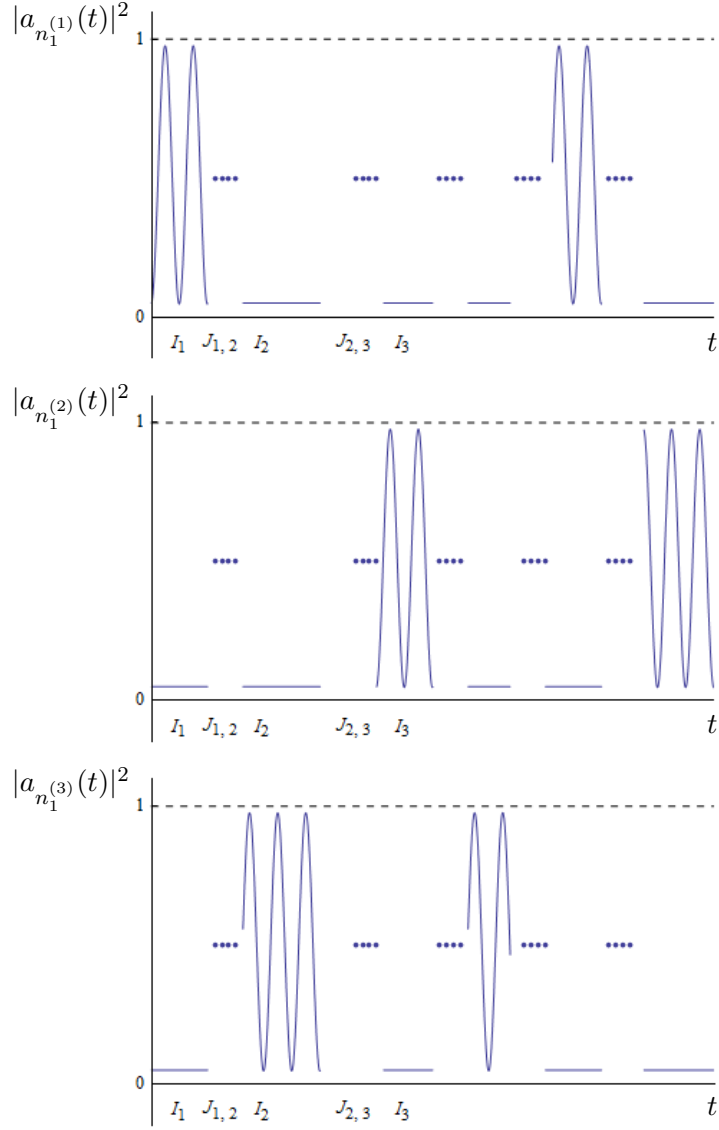
\begin{figure}[b!] 
\begin{overpic}[scale=0.9]{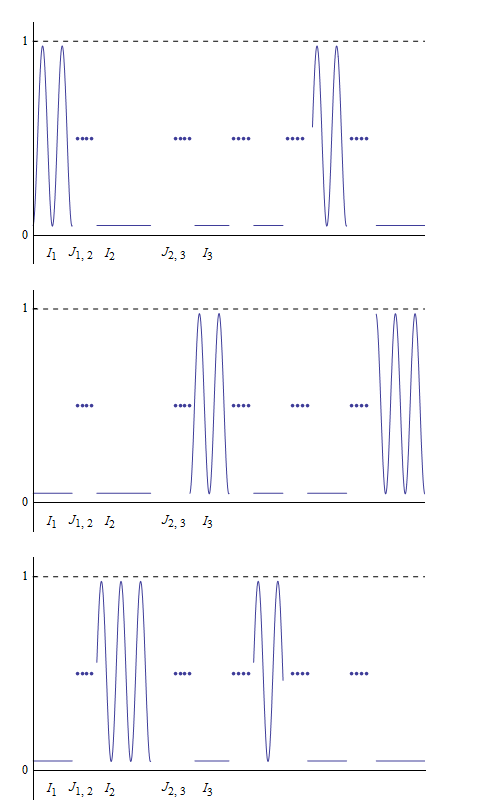}
\put(-8,96){$|a_{n_{1}^{(1)}}(t)|^2$ }
\put(52,68){$t$}

\put(-8,63){$|a_{n_{1}^{(2)}}(t)|^2$ }
\put(52,35){$t$}

\put(-8,30){$|a_{n_{1}^{(3)}}(t)|^2$ }
\put(52,2){$t$}
\end{overpic}
\caption{An example of the evolution of $|a_{n_1}^{(j)}(t)|^2$, $j=1, 2, 3$  of a solution obtained by Theorem \ref{thm:traveling_periodic_beating} as time evolves. We consider $N=3$ and the sequence of modes which are ``activated'' is $\omega=\{ 1, 3, 2, 3, 1, 2\ldots\}$.}
\label{fig:travel}
\end{figure}

The orbits of Theorem \ref{thm:geometric} are obtained by shadowing a sequence
of invariant tori (periodic orbits for a suitable symplectic reduction)
connected by transverse heteroclinic orbits. Then, the
\emph{beating-time} intervals are the  intervals where the orbit is in a
small neighborhood of each of the periodic orbits. In this regime, (the moduli
of) some modes oscillate periodically, whereas the others are at rest.  The
\emph{transition-time} intervals correspond to  intervals in which the orbit
is ``traveling'' along a heteroclinic orbit and is ``far'' from all periodic
orbits. In this regime, all modes undergo a drastic change (see Figure \ref{fig:travel}).
\subsection{Some comments on Theorems \ref{TeoWaveBeam} and
\ref{thm:traveling_periodic_beating}}

\smallskip

\noindent\emph{Hartree equation.} Similar results hold true also for  the
Hartree equation
\begin{equation}\label{Hartree}
\mathrm{i} u_t=\Delta u+(V\star \lvert u \rvert^2)\,u,\qquad u=u(t, x), \quad
t\in\mathbb{R}, \quad x\in\mathbb{T}^2
\end{equation}
with a convolution potential $V(x)=\sum_{j\in\mathbb{Z}^2} V_j\,e^{\mathrm{i} j
\cdot x}$ such that
\begin{equation}\label{def:potential}
V\colon \mathbb{T}^2\to\mathbb{R}, \quad V(x)=V(-x)
\end{equation}
and assuming the following hypothesis. Once fixed the set
$\Lambda\subset\mathbb{Z}^2$ of resonant tuples, the Fourier coefficients $V_j$
of the  potential
with $j=n_1-n_2$ for some $n_1, n_2\in\Lambda$ satisfy
\begin{equation}\label{cond:potential}
\quad V_j=1+\eps\gamma_j \quad \text{with}\quad \eps\ll 1.
\end{equation}
Assume that  the coefficients $\gamma_j$ satisfy a codimension 1 non-degeneracy condition
and take $\eps$ small enough. Then, the Hartree
equation has solutions of the form
\[
 u(t, x)=\de\sum_{n\in \Lambda} a_{n}(t)\,e^{\mathrm{i} n \cdot x}+R(t, x)
\]
where the first order $\{a_{n}\}$ and the remainder $R$ satisfy the statements
given either in  Theorem \ref{TeoWaveBeam} (where $R\rightsquigarrow R_3$) or
\ref{thm:traveling_periodic_beating} (where $R\rightsquigarrow R_4$) .

\smallskip
\noindent\emph{Smale Horseshoes in PDEs:} Theorem \ref{thm:geometric} provides a
Smale Horseshoe for the Birkhoff normal form. This invariant set is partially
hyperbolic and partially elliptic if considered in the whole infinite
dimensional phase space. This is the reason why, a priori, this invariant set is
not persistent for the full equations \eqref{Wave}, \eqref{Beam},
\eqref{Hartree}. As
far as the authors know, the existence of Smale horseshoes in Hamiltonian PDEs
has been mostly obtained by adding dissipation to the equation which make these
sets become fully hyperbolic
\cite{HolmesM81,Li99,BertiC02,BatelliF05,BatelliG08}.

%

\smallskip

\noindent\emph{Non-integrability of $\mathcal{N}$ on $V_{\Lambda}$:}
 Theorem \ref{thm:geometric} (and therefore Theorems \ref{TeoWaveBeam} and
\ref{thm:traveling_periodic_beating}) relies on the fact that the restriction of
$\mathcal{N}$ on $V_{\Lambda}$ is not integrable and admits invariant tori with
transverse homoclinic orbits which lead to chaotic behavior. If one tries to mimic the proof of Theorems \ref{TeoWaveBeam} and
\ref{thm:traveling_periodic_beating} for the cubic NLS
\begin{equation*}
 i u_t=\Delta u-|u|^2 u,\quad x\in \mathbb{T}^2,
\end{equation*}
one fails to obtain chaotic dynamics for the normal form $\mathcal{N}$ on $V_{\Lambda}$. This is due to the fact that, for the sets $\Lambda$ considered in Theorems \ref{TeoWaveBeam} and \ref{thm:traveling_periodic_beating}, the dynamics of $\mathcal{N}_{|_{V_{\Lambda}}}$ is integrable.
Indeed the dynamics of the normal form restricted to each resonant tuple is decoupled and integrable for NLS (recall that $\Lambda$ is chosen as union of disjoint resonat tuples). The integrability prevents the existence of homoclinic or heteroclinic transverse intersections and chaotic motions. \\
Certainly we expect the normal form of NLS to be non-integrable, but one should consider solutions supported on sets $\Lambda$ different from the ones considered in Theorems \ref{TeoWaveBeam} and
\ref{thm:traveling_periodic_beating} in order to capture the chaotic behavior.

\smallskip

\noindent\emph{Defocusing and Focusing equations:} To simplify the exposition,
the theorems above only refer to the defocusing equations \eqref{Wave} and
\eqref{Beam}. However, it can be checked that the sign of the nonlinearity does
not play any role and therefore, Theorems \ref{thm:geometric}, \ref{TeoWaveBeam}
and \ref{thm:traveling_periodic_beating} also apply to the focusing equations
 \[
  u_{tt}-\Delta u-u^3=0,\qquad  u_{tt}+\Delta^2 u-u^3=0.
 \]

\smallskip

\noindent\emph{Transfer of energy and growth of Sobolev norms:}
The solutions of the Wave equation \eqref{Wave}/Beam equation
\eqref{Beam}/Hartree equation \eqref{Hartree} obtained in Theorem
\ref{thm:traveling_periodic_beating} undergo certain transfer of energy between
modes. Unfortunately, such transfer of energy do not lead to growth of Sobolev
norms \cite{Bourgain00b, CKSTT, GuardiaK12}. In \cite{CKSTT}, the authors obtain
orbits undergoing growth of Sobolev norms for the defocusing NLS equation on
$\mathbb{T}^2$. One of the key points of their proof is to construct, for the
Birkhoff normal form truncation, a  chain of invariant tori (periodic orbits in
certain symplectic reduction, named \emph{toy model}) which are connected by
\emph{non-transverse} heteroclinic orbits (see also \cite{DelshamsSZ18} for a thorough analysis of non-transverse shadowing arguments). In order to obtain such connections, they rely on the fact that the toy model is integrable once restricted to certain invariant subspace (called \emph{two generations model} in \cite{CKSTT}). Then the orbits undergoing growth of Sobolev norms are well approximated by orbits which shadow (follow closely) this chain of periodic orbits. 
This approach does not work with the Wave \eqref{Wave}, Beam  \eqref{Beam} and Hartree \eqref{Hartree} equations. First of all, in these equations the two generations model is not integrable (for the Hartree equation it is not for a generic potential), and the models we consider are carefully chosen so that they are close to integrable, and therefore can be analyzed through perturbative methods; for the Wave and Beam equation, to be close to integrable we have to choose the modes in $\Lambda$ with very similar modulus, and it seems difficult to use the analysis done in this paper to construct orbits undergoing growth of Sobolev norms. For the Hartree equation, one should expect that the ideas developed in this paper could lead to growth of Sobolev norms for a generic potential satisfying \eqref{def:potential}, \eqref{cond:potential}. A second fundamental difference between NLS and the PDEs considered in this paper is that the chain of tori connected by  heteroclinic connections considered in \cite{CKSTT} is \emph{not structurally stable} since to have such heteroclinic connections one certainly needs that the connected invariant tori belong to the same level of energy; indeed, for the Hartree equation \eqref{Hartree} with a
generic potential $V$ the tori considered in \cite{CKSTT} belong to different
level of energy and the same happens for the Wave and Beam equations for a
generic choice of resonant tuples. The tori considered in  Theorem
\ref{thm:geometric} are radically different from those in \cite{CKSTT}: these
tori come in families of higher dimension which are transverse to the first
integrals, and are connected by transverse heteroclinic orbits. We believe that
such objects could play a role if one wants to implement \cite{CKSTT} to other
PDEs.

\section*{Acknowledgments}

This project has received funding from the European Research Council (ERC) under
the European Union's Horizon 2020 research and innovation programme (grant
agreement No 757802). P.M. has been partially funded by the Spanish Government
MINECO-FEDER grant PGC2018-100928-B-I00. M. G and S. P. have been also partially
supported by
the Spanish MINECO-FEDER Grant PGC2018-098676-B-100 (AEI/FEDER/UE) and by the
Catalan Institution for Research
and Advanced Studies via an ICREA Academia Prize 2019.  F.G., M. G, P. M. and
S.P. have  been also partially supported by the Catalan grant 2017SGR1049. S.P.
acknowledges financial support from the Spanish ``Ministerio de Ciencia,
Innovaci\'on y
Universidades'', through the Mar\'{\i}a de Maeztu Programme for Units of
Excellence (2015-
2019) and the Barcelona Graduate School of Mathematics.

\section{Strategy of the proofs} \label{sec:Proof}

The general argument we use in the proofs of Theorems~\ref{TeoWaveBeam}
and~\ref{thm:traveling_periodic_beating} follows some of the ideas in the
literature~\cite{CKSTT,GuardiaK12,Guardia14,HausProcesi,GuardiaHP16}. The steps
of the proof of both theorems are the following.

\smallskip

\noindent\textbf{Step 1: }The PDEs under consideration have a Hamiltonian
structure. Let us denote by $H$  the Hamiltonian, which can be written as
\[
 H=H^{(2)}+H^{(4)},
\]
where $H^{(i)}$ are homogenous polynomials of degree $i$.
Given a finite subset $\Lambda\subset\mathbb{Z}^2$ of resonant modes, to be
chosen later, we apply a \emph{weak normal form} scheme to $H$.
 Let us briefly describe such procedure. Consider the splitting $u=v+z$ where
$v$ is the projection of $u$ on the finite dimensional subspace $V_{\Lambda}$ in
\eqref{def:Vlambda}
 and $z$ is the projection on the $L^2$-orthogonal of this subspace.
Let us denote by
$H^{(4, k)}$, $0\le k\le 4$, the projection of $H^{(4)}$ onto the monomials of
the form $v^{4-k}\,z^k$. We construct a change of coordinates $\Gamma$ which
normalizes all the Hamiltonian terms of degree $4$ which are independent or
linear in the normal variables, namely $H^{(4, 0)}$ and $H^{(4, 1)}$. This
normal form procedure is enough for our purposes. The terms in $H^{(4, 0)}$ and $H^{(4, 1)}$ that cannot be removed or normalized (the resonant terms) are the ones Fourier supported on the following resonances
\[
\begin{cases}
\sum_{i=1}^4 \sigma_i j_i=0,\\
\sum_{i=1}^4 \sigma_i \omega(j_i)=0,
\end{cases}
\]
where $\sigma_i\in \{ \pm 1 \}$, $j_1,\dots, j_4\in\mathbb{Z}^2_{\odd}$ and at
most one of them is outside $\Lambda$, $\omega(j_i)$ is the linear frequency of
oscillation of the $j_i$-th mode (see \eqref{def:variational}).

 Thanks to the momentum conservation, expressed by the first equation of the
above system, the terms $H^{(4, 0)}$ and $H^{(4, 1)}$ are given by a finite sum
of monomials.  As a consequence, the normalizing transformation $\Gamma$ is
well
defined as the time-one flow map of an ODE. We observe that  expressions of the
form $\sum_{i=1}^4 \sigma_i \omega(j_i)$ appear at the denominator of the
Fourier coefficients of the Hamiltonian that generates the Birkhoff map. These
quantities in general may accumulate to zero, as in the case of the Wave
equation, but by the finiteness argument they are bounded from below by a
constant that depends on $\Lambda$.

We choose $\Lambda$ such that there are no resonant monomials supported on just
one mode outside $\Lambda$, i.e. $H^{(4, 1)}=0$. This means that the resonant
Hamiltonian of degree
$4$, we call it $H_{\mathrm{Birk}}$, is given by the normalization of $H^{(4,
0)}$. Then,
\[
H\circ \Gamma=H^{(2)}+H_{\mathrm{Birk}}+R,
\]
where $R$ is a remainder which contains degree 4 monomials supported at least on
two normal modes or monomials with degree greater or equal than $6$. Such terms
shall be considered as a small perturbation of the truncated normal form
$H^{(2)}+H_{\mathrm{Birk}}$. Moreover it is easy to see that $V_{\Lambda}$
 is invariant by the flow of the truncated normal form. The \emph{resonant
model} is obtained by considering the restriction of the normal form on
$V_{\Lambda}$. In the following step we discuss how we construct particular
orbits for this system.\\

\noindent\textbf{Step 2: } This step is the core of the paper and can be divided
as follows.

\noindent\emph{Construction of the set $\Lambda$:} The set $\Lambda$ is chosen
in such a way that its associated subspace of modes $V_{\Lambda}$ is invariant
by the flow of the resonant model, but of course satisfies other requirements.
Its precise definition depends on the PDE model we consider, but all three
instances (Wave, Beam and Hartree equations) of the set $\Lambda$ share some
common features. They have exactly $4N$ elements, $N\geq 2$, which, using the
terminology introduced in~\cite{CKSTT}, encompass two generations. The elements
of the set $\Lambda$ are organized in groups of four,
 pairwise disjoint, each of them forming a ``resonant'' parallelogram. The
choice of the modes is such that each individual parallelogram
is invariant; moreover, the dynamics of a single parallelogram is integrable.

\noindent\emph{The dynamics of the finite dimensional model:}
 We choose the modes in $\Lambda$ in such a way that the dynamics of the
resonant model is close to integrable, where closeness to integrability is
measured through some parameter $\varepsilon$ (see Remark \ref{rmk:Lambda}). In
the case of the Wave and Beam equation the case $\varepsilon=0$ corresponds to
consider all the modes in $\Lambda\subset\mathbb{Z}^2$ on a circle. We obtain
the nearly integrability by placing the resonant modes close to that circle.

After a symplectic reduction the Hamiltonian of the resonant model has
$N$-degrees of freedom and it reads

\begin{equation}\label{redhamN}
\begin{aligned}
\mathcal{H}(\psi_i, K_i) &=\mathcal{H}_0(\psi_i, K_i)+\varepsilon
\mathcal{H}_1(\psi_i, K_i), \quad i=1, \dots, N,
\end{aligned}
\end{equation}
with
\[
 \begin{split}
\mathcal{H}_0 &= \sum_{j=1}^N K_j (1-K_j)(1 +2  \,\cos(\psi_j))\\
\mathcal{H}_1&=\sum_{j=1}^N \left(a_{j} K_j + b_{j} K^2_j+c_j K_j(1-K_j)
\cos\psi_j\right) + \sum_{i, j=1, i< j}^N d_{i j} K_iK_j,
 \end{split}
\]
for some coefficients $a_j, b_j, c_j, d_{i j}\in\mathbb{R}$ depending on the
choice of $\Lambda$ (and on the potential $V$ in the case of the Hartree
equation).

 The unperturbed system (where $\varepsilon = 0$) possesses certain invariant objects, namely  hyperbolic fixed points and hyperbolic periodic orbits, whose invariant manifolds form heteroclinic separatrices.
   One can find a similar picture by considering the NLS Hamiltonian restricted to the subspace generated by $4$ resonant modes which form a non-degenerate rectangle (or many copies of such system).

The chaotic dynamics that build originates from the transverse
homoclinic connections that arise for $0<\eps\ll 1$. Note that these
homoclinics \emph{do not exist} in the unperturbed case $\eps=0$ (where the
invariant manifolds form heteroclinic connections instead of homoclinic). To
prove the existence of such homoclinics one cannot apply directly the standard
Melnikov theory \cite{Melnikov63}, but suitable generalizations of it.
With similar arguments, we construct
heteroclinic connections between objects which are not connected in the
unperturbed case.

The conditions imposed on the  set $\Lambda$
(namely, suitable Melnikov-like functions having non-degenerate zeros)
ensures that these transverse  homoclinic and heteroclinic connections exist
for small $\varepsilon \neq 0$.\\


\noindent\emph{The infinite symbols Smale horseshoe:}  The orbits in
Theorem~\ref{TeoWaveBeam} give rise from a horseshoe of infinite symbols that
can be constructed close to a hyperbolic periodic orbit of \eqref{redhamN}
 whose invariant manifolds intersect transversally. The construction of this
horseshoe follows the ideas in~\cite{Moser01}. In this horseshoe, each symbol
encodes the time in which the trajectory comes close to the
periodic orbit.
The horseshoe and its dynamics can be described as follows. Let
$\Gamma=\{1,2,3, \dots\}$ be a denumerable set of symbols and
\[
\Sigma = \{s= (\dots, s_1, s_0, s_1, \dots) \mid s_i \in \Gamma, \; i \in \N\},
\]
the space of bi-infinite \emph{sequences}, with the product topology. Notice that, unlike what happens when $\Gamma$ is a finite set, $\Sigma$ is not compact. The shift $\sigma: \Sigma \to \Sigma$ is the homeomorphism
on $\Sigma$ defined by $(\sigma(s))_{i} = s_{i-1}$.\\

\noindent \emph{Shadowing of a sequence of periodic orbits:} The orbits in
Theorem~\ref{thm:traveling_periodic_beating}
travel along a chain of periodic orbits connected by transverse heteroclinic orbits, following the diffusion mechanism
described originally by Arnold~\cite{Arnold64}. This mechanism consists of a
sequence - finite or infinite - of partially hyperbolic periodic orbits,
$\{\TT_i\}_{i\in I}$, $I\subset\mathbb{N}$, such that the unstable manifold of
$\TT_i$, $W^u(\TT_i)$, intersects \emph{transversally} the stable manifold of
$\TT_{i+1}$, $W^s(\TT_{i+1})$. Here, since the system we are considering is
autonomous, \emph{transversally} means transversality in the energy level, which
 implies that the intersection of the manifolds is, locally, a single
heteroclinic orbit. If a nondegeneracy condition is met, this transversality is
sufficient to apply a Lambda Lemma that implies that
$W^u(\TT_{i+1}) \subset \overline{W^u(\TT_{i})}$ (see~\cite{FontichM01}), which in turn implies that for any $i,j \in I$,
$i< j$, $W^u(\TT_{j}) \subset \overline{W^u(\TT_{i})}$. One can then choose arbitrary small neighborhoods of the tori $\TT_i$ and orbits that visit these neighborhoods according to an increasing sequence of times.

It is worth to remark that the orbits found in the resonant model do exist for
any positive time.  In the case of the horseshoe with infinite symbols, one
obtains orbits that arrive at randomly chosen times to a neighborhood of the
periodic orbit. In the case  of the diffusion orbits,
one obtain solutions that wander along the chain of periodic orbits for any
positive time, and can be chosen to arrive closer and closer to each periodic
orbit.

\smallskip

\noindent\textbf{Step 3:} The last step of the proof consists on finding a true solution of each PDE shadowing for long enough time the chosen solution of the resonant model. This is accomplished by a standard Gronwall and bootstrap argument.

\bibliography{references}
\bibliographystyle{siam}


\end{document}